\documentclass[11pt]{amsart}

\usepackage{amsmath}
\usepackage{amssymb}
\usepackage{mathrsfs}
\usepackage{comment}
\usepackage{color}
\usepackage[colorlinks,citecolor=blue,urlcolor=black,linkcolor=black]{hyperref}

\newtheorem{theorem}{Theorem}[section]
\newtheorem{claim}[theorem]{Claim}

\newtheorem{lemma}[theorem]{Lemma}

\newtheorem{proposition}[theorem]{Proposition}

\theoremstyle{definition}

\newtheorem{question}[theorem]{Question}

\theoremstyle{remark}
\newtheorem{remark}[theorem]{Remark}

\newcount\skewfactor
\def\mathunderaccent#1#2 {\let\theaccent#1\skewfactor#2
\mathpalette\putaccentunder}
\def\putaccentunder#1#2{\oalign{$#1#2$\crcr\hidewidth
\vbox to.2ex{\hbox{$#1\skew\skewfactor\theaccent{}$}\vss}\hidewidth}}
\def\name{\mathunderaccent\tilde-3 }

\def\smallbox#1{\leavevmode\thinspace\hbox{\vrule\vtop{\vbox
   {\hrule\kern1pt\hbox{\vphantom{\tt/}\thinspace{\tt#1}\thinspace}}
   \kern1pt\hrule}\vrule}\thinspace}

\newcommand{\cf}{{\rm cf}}

\def\qedref#1{$\qed_{\reforiginal{#1}}$}

\title{Hungarian cubes}
\author{Shimon Garti}
\address{Einstein Institute of Mathematics,
 The Hebrew University of Jerusalem,
 Jerusalem 91904, Israel}
\email{shimon.garty@mail.huji.ac.il}

\thanks{Research supported by ISF grant no. 2320/23}

\begin{document}
\let\labeloriginal\label
\let\reforiginal\ref
\def\ref#1{\reforiginal{#1}}
\def\label#1{\labeloriginal{#1}}

\begin{abstract}
We prove the consistency of the relation $\left( \begin{smallmatrix} \nu \\ \mu \\ \lambda \end{smallmatrix} \right) \rightarrow \left( \begin{smallmatrix} \nu \\ \mu \\ \lambda \end{smallmatrix} \right)$ when $\lambda<\mu=\cf(\mu)<\nu=\cf(\nu)\leq 2^\mu$.
\end{abstract}

\maketitle

\newpage

\section{Introduction}

The polarized cube relation $\left( \begin{smallmatrix} \alpha \\ \beta \\ \gamma \end{smallmatrix} \right) \rightarrow \left( \begin{smallmatrix} \varepsilon \\ \zeta \\ \eta \end{smallmatrix} \right)$ says that for every coloring $d:\alpha\times\beta\times\gamma\rightarrow 2$ one can find $A\subseteq\alpha,B\subseteq\beta$ and $C\subseteq\gamma$ such that ${\rm otp}(A)=\varepsilon,{\rm otp}(B)=\zeta,{\rm otp}(C)=\eta$ and $d\upharpoonright(A\times B\times C)$ is constant.
In his famous remark from \cite{pdf}, Pierre de Fermat says that \emph{Cubum autem in duos cubos}, that is: a cube cannot split into two cubes.
Though true for natural numbers, as proved by Wiles in \cite{MR1333035}, it may not hold at infinite cardinals.
Our goal is to prove positive combinatorial relations for cubes whose edges are infinite cardinals.

A necessary condition for a positive cube relation is the positive standard relation at all three pairs of cardinals mentioned in the cube.
Thus if $\binom{\alpha}{\beta}\nrightarrow\binom{\varepsilon}{\zeta}$ or $\binom{\beta}{\gamma}\nrightarrow\binom{\zeta}{\eta}$ or $\binom{\alpha}{\gamma}\nrightarrow\binom{\varepsilon}{\eta}$ then $\left( \begin{smallmatrix} \alpha \\ \beta \\ \gamma \end{smallmatrix} \right) \nrightarrow \left( \begin{smallmatrix} \varepsilon \\ \zeta \\ \eta \end{smallmatrix} \right)$.
Moreover, positive relations at these pairs are necessary but insufficient and usually something very strong is needed for a positive cube relation in non-trivial cases.

The case $\left( \begin{smallmatrix} \alpha \\ \beta \\ \gamma \end{smallmatrix} \right) \rightarrow \left( \begin{smallmatrix} \alpha \\ \beta \\ \gamma \end{smallmatrix} \right)$ in which the required size of the monochromatic cube equals that of the domain of the coloring is called \emph{the strong cube relation}.
Since $\binom{\kappa}{\kappa}\nrightarrow\binom{\kappa}{\kappa}$ we see that the strong relation $\left( \begin{smallmatrix} \nu \\ \mu \\ \lambda \end{smallmatrix} \right) \rightarrow \left( \begin{smallmatrix} \nu \\ \mu \\ \lambda \end{smallmatrix} \right)$ is possible only if $\lambda,\mu,\nu$ are all distinct.
Our convention is that $\lambda\leq\mu\leq\nu$ so we may assume from now on that $\lambda<\mu<\nu$.
If $\nu=\cf(\nu)>2^\mu$ then $\nu$ adds nothing to the validity of the cube relation and $\left( \begin{smallmatrix} \nu \\ \mu \\ \lambda \end{smallmatrix} \right) \rightarrow \left( \begin{smallmatrix} \nu \\ \mu \\ \lambda \end{smallmatrix} \right)$ iff $\binom{\mu}{\lambda}\rightarrow\binom{\mu}{\lambda}$ thus the cube problem reduces to a standard polarized relation.

The opposite situation in which $2^\lambda$ is relatively small does not reduce immediately to the case $\binom{\nu}{\mu}\rightarrow\binom{\nu}{\mu}$.
An easy case is $2^\lambda<\kappa<\mu<\nu$ where $\kappa$ is strongly compact and $\binom{\nu}{\mu}\rightarrow\binom{\nu}{\mu}$ can be witnessed by sets which belong to some prescribed $\lambda$-complete ultrafilters over $\mu$ and $\nu$.
In this case one has $\left( \begin{smallmatrix} \nu \\ \mu \\ \lambda \end{smallmatrix} \right) \rightarrow \left( \begin{smallmatrix} \nu \\ \mu \\ \lambda \end{smallmatrix} \right)$ but without this assumption the cube relation is more involved even if $2^\lambda$ is small.
It seems that one has to distinguish between the case of $2^\lambda<\mu$ and the case of $2^\lambda\geq\mu$, and in both cases the cube relation $\left( \begin{smallmatrix} \nu \\ \mu \\ \lambda \end{smallmatrix} \right) \rightarrow \left( \begin{smallmatrix} \nu \\ \mu \\ \lambda \end{smallmatrix} \right)$ seems to be interesting.

The most challenging case seems to be $\nu\leq 2^\lambda$, see Question \ref{qbasa}.
In this paper we obtain the strong cube relation $\left( \begin{smallmatrix} \nu \\ \mu \\ \lambda \end{smallmatrix} \right) \rightarrow \left( \begin{smallmatrix} \nu \\ \mu \\ \lambda \end{smallmatrix} \right)$ where $\nu\leq 2^\mu$ but $2^\lambda$ is relatively small.
Thus the main result is the consistency of $\left( \begin{smallmatrix} \nu \\ \mu \\ \lambda \end{smallmatrix} \right) \rightarrow \left( \begin{smallmatrix} \nu \\ \mu \\ \lambda \end{smallmatrix} \right)$ when $\lambda<\mu=\cf(\mu)<\nu=\cf(\nu)=2^\mu$ in Theorem \ref{thmforcemt}.
We indicate, however, that the combinatorial theorem that we prove is phrased in a general way and might be applied to cases in which $\nu\leq 2^\lambda$ provided that the pertinent assumptions are forceable.

We consider in this paper only strong relations, namely relations of the form $\left( \begin{smallmatrix} \alpha \\ \beta \\ \gamma \end{smallmatrix} \right) \rightarrow \left( \begin{smallmatrix} \alpha \\ \beta \\ \gamma \end{smallmatrix} \right)$.
The literature concerning cube polarized relations is sparse.
An old negative relation under \textsf{GCH} appeared in \cite{MR0048517}.
A famous problem from \cite{MR0280381} is whether $\left( \begin{smallmatrix} \aleph_1 \\ \aleph_1 \\ \aleph_1 \end{smallmatrix} \right) \rightarrow \left( \begin{smallmatrix} \aleph_0 \\ \aleph_0 \\ \aleph_0 \end{smallmatrix} \right)$ and more generally whether $\left( \begin{smallmatrix} \kappa^+ \\ \kappa^+ \\ \kappa^+ \end{smallmatrix} \right) \rightarrow \left( \begin{smallmatrix} \kappa \\ \kappa \\ \kappa \end{smallmatrix} \right)$ for some $\kappa$, see \cite[p. 110]{MR3075383}.
For interesting results concerning cubes of this type we refer to \cite{MR1297180}.
However, cube relations in which the size of the monochromatic product is smaller than the size of the domain will not be discussed in the current paper.

Our notation is mostly standard, and follows \cite{MR795592}.
We shall use the Jerusalem forcing notation, thus $p\leq{q}$ means that $p$ is weaker than $q$.
For information about indecomposable filters we refer to \cite{MR727789} and \cite{MR865928}.
We suggest \cite{MR2768693} and \cite{MR1318912} for background in pcf theory and \cite{MR3075383} for basic background in partition calculus.
Many thanks to the referee of the paper for a lot of helpful suggestions which improved the accuracy and the readability of the article.

\newpage 

\section{Preliminaries}

In this section we state several results, to be used in the next sections.
We commence with a statement about the ultrafilter number $\mathfrak{u}_\mu$, the minimal size of a base for some uniform ultrafilter over $\mu$.
Recall that a filter $\mathscr{F}$ over $\kappa$ is $\theta$-indecomposable iff for every partition $(A_\alpha:\alpha\in\theta)$ of $\kappa$ one can find an index set $I\in[\theta]^{<\theta}$ so that $\bigcup_{\alpha\in I}A_\alpha\in\mathscr{F}$.

\begin{proposition}
  \label{propraghavanshelah} Assume that:
  \begin{enumerate}
    \item [$(a)$] $\kappa$ is supercompact and Laver-indestructible.
    \item [$(b)$] $\kappa^+<\theta=\cf(\theta)<\mu=\cf(\mu)$.
    \item [$(c)$] $\mu'$ is a strong limit singular cardinal.
    \item [$(d)$] $\theta=\cf(\mu')$ and $\mu<\mu'$.
    \item [$(e)$] there exists a $\theta$-indecomposable filter over $\mu$.
    \item [$(f)$] $\nu=\cf(\nu)>\mu'$.
  \end{enumerate}
  Then $\mathbb{C}_{\theta\nu}$ forces $\mathfrak{u}_\mu<2^\mu=\nu$ while preserving the supercompactness of $\kappa$.
\end{proposition}

\par\noindent\emph{Proof}. \newline
By \cite[Theorem 7]{MR4081063}, where $\mu'$ here stands for $\mu$ there.
Observe that $\kappa$ remains supercompact since $\mathbb{C}_{\theta\nu}$ is $\kappa$-directed-closed.

\hfill \qedref{propraghavanshelah}

We turn now to a short discussion concerning Mathias forcing from \cite{MR491197} and the generalized Mathias forcing from \cite[Definition 3.1]{MR2927607}.
Let $\mathscr{U}$ be a non-principal ultrafilter over $\omega$.
Mathias forcing relativized to $\mathscr{U}$ is denoted by $\mathbb{M}_{\mathscr{U}}$.
A condition $p\in\mathbb{M}_{\mathscr{U}}$ is a pair $(s^p,A^p)$ where $s^p\in[\omega]^{<\omega}, A^p\in\mathscr{U}$ and $\max(s^p)<\min(A^p)$.
If $p,q\in\mathbb{M}_{\mathscr{U}}$ then $p\leq q$ iff $s^p\subseteq s^q, A^p\supseteq A^q$ and $s^q-s^p\subseteq A^p$.
If $G\subseteq\mathbb{M}_{\mathscr{U}}$ is $V$-generic then $\bigcup\{s^p:p\in G\}$ is a subset of $\omega$ called \emph{a Mathias real}.

The advantage of relativizing Mathias forcing to an ultrafilter is that it becomes $ccc$ and hence can be iterated with a finite support without collapsing cardinals.
An iteration of $\mathbb{M}_{\mathscr{U}}$ over a model of \textsf{CH} may change the splitting number $\mathfrak{s}$ as described in the following proposition.

\begin{claim}
  \label{clmsplitting} Assume that:
  \begin{enumerate}
    \item [$(a)$] $\mathsf{CH}$ holds in the ground model.
    \item [$(b)$] $\theta=\cf(\Upsilon)>\omega$.
    \item [$(c)$] $\mathbb{P}$ is a finite support iteration of Mathias forcing notions relativized to some ultrafilter.
    \item [$(d)$] The length of the iteration is $\Upsilon$.
    \item [$(e)$] $G\subseteq\mathbb{P}$ is $V$-generic.
  \end{enumerate}
  Then $\mathfrak{s}=\theta$ in $V[G]$.
\end{claim}

\par\noindent\emph{Proof}. \newline
For every $\alpha\in\Upsilon$ let $S_\alpha=\{x^\alpha_n:n\in\omega\}$ be the $\alpha$th Mathias real, and let $S^0_\alpha=\{x^\alpha_{2n}:n\in\omega\}$.
Let $(\alpha_i:i\in\theta)$ be an increasing cofinal sequence of ordinals in $\Upsilon$.
Let $\mathcal{S}=\{S^0_{\alpha_i}:i\in\theta\}$, so $|\mathcal{S}|\leq\theta$.
As a first step, we show that $\mathcal{S}$ is a splitting family in $V[G]$, thus $V[G]\models\mathfrak{s}\leq\theta$.
We will show, then, that any smaller family is not a splitting family in $V[G]$, thus $V[G]\models\mathfrak{s}\geq\theta$.
Together, the statement will be proved.

By the properties of Mathias forcing, $(S_{\alpha_i}:i\in\theta)$ is $\subseteq^*$-decreasing.
Let $\mathscr{W}$ be the filter generated by $(S_{\alpha_i}:i\in\theta)$.
We claim that $\mathscr{W}$ is actually an ultrafilter.
Indeed, if $A\in V[G]\cap[\omega]^\omega$ then $A$ belong to some proper initial segment of the iteration, since $\theta=\cf(\theta)>|A|$ and $\mathbb{P}$ is a finite support iteration of $ccc$ forcing notions.
Hence for some $i\in\theta$ one has either $S_{\alpha_i}\subseteq^*{A}$ or $S_{\alpha_i}\subseteq^*{\omega-A}$.
This means that either $A\in\mathscr{W}$ or $\omega-A\in\mathscr{W}$, so $\mathscr{W}$ is an ultrafilter over $\omega$ in $V[G]$.
Fix $B\in V[G]\cap[\omega]^\omega$.
Without loss of generality, $B\in\mathscr{W}$ and hence $S_{\alpha_i}\subseteq^*{B}$ for some $i\in\theta$.
It follows that $S^0_{\alpha_i}$ splits $B$, thus $\mathcal{S}$ is a splitting family as desired.

Now suppose that $\theta'<\theta$ and $\mathcal{T}=\{T_\gamma:\gamma\in\theta'\}\subseteq[\omega]^\omega$.
We may assume that $\mathcal{T}\subseteq\mathscr{W}$ (indeed, let $\mathcal{T}'=\{T'_\gamma:\gamma\in\theta'\}$ where $T'_\gamma=T_\gamma$ if $T_\gamma\in\mathscr{W}$ and $T'_\gamma=\omega-T_\gamma$ if $T_\gamma\notin\mathscr{W}$.
Now $\mathcal{T}'\subseteq\mathscr{W}$ and $\mathcal{T}$ is a splitting family iff $\mathcal{T}'$ is a splitting family).
For every $\gamma\in\theta'$ there exists $i_\gamma\in\theta$ so that $S_{\alpha_{i_\gamma}}\subseteq^* T_\gamma$.
Since $\theta'<\theta=\cf(\Upsilon)$ there is $i\in\theta$ so that $\alpha_i>\bigcup_{\gamma\in\theta'}\alpha_{i_\gamma}$.
Consequently, $S_{\alpha_i}\subseteq^* T_\gamma$ for every $\gamma\in\theta'$, thus $S_{\alpha_i}$ is not split by $\mathcal{T}$.
We conclude that $\mathfrak{s}\geq\theta$ in $V[G]$, so we are done.

\hfill \qedref{clmsplitting}

It follows from the above claim that in $V[G]$ if $\omega<\cf(\mu)\leq\mu<\mathfrak{s}$ then $\binom{\mu}{\omega}\rightarrow\binom{\mu}{\omega}$.
But this forcing notion has a similar effect on cardinals above $\mathfrak{s}$.
If one begins the iteration with blowing up $2^\omega$ to $\nu$ and then iterates $\mathbb{M}_{\mathscr{U}}$ say $\lambda$-many times where $\omega<\lambda=\cf(\lambda)<\nu$ then $\binom{\mu}{\omega}\rightarrow\binom{\mu}{\omega}$ holds in $V[G]$ for every $\lambda<\cf(\mu)\leq\mu\leq\nu$.
The reason is that the $\lambda$-sequence of Mathias reals forms a $\subseteq^*$-decreasing sequence of subsets of $\omega$ and thus generates the non-principal ultrafilter $\mathscr{W}$ over $\omega$ whose character is smaller than $\mu$.
In particular, if $2^\omega=\nu$ and $\cf(\nu)>\lambda$ then $\binom{\nu}{\omega}\rightarrow\binom{\nu}{\omega}$ holds in $V[G]$.
This will be used in Claim \ref{clmufnumber}.

Mathias forcing generalizes to uncountable cardinals in lieu of $\omega$.
Suppose that $\lambda$ is measurable and $\mathscr{U}$ is a $\lambda$-complete ultrafilter over $\lambda$.
A single step of the generalized Mathias forcing consists of conditions of the form $(s^p,A^p)$ where $s^p\in[\lambda]^{<\lambda}, A^p\in\mathscr{U}$ and $\sup(s^p)<\min(A^p)$.
The $\lambda$-completeness of $\mathscr{U}$ is needed for $\mathbb{M}_{\mathscr{U}}$ to be $\lambda$-closed, hence measurability is indispensable.
In order to iterate $\mathbb{M}_{\mathscr{U}}$ one has to ensure that $\lambda$ remains measurable at each step of the iteration.
If one begins with a Laver indestructible supercompact cardinal $\lambda$ then this can be done since $\mathbb{M}_{\mathscr{U}}$ will be $\lambda$-directed-closed.
If one forces $2^\lambda=\nu=\cf(\nu)$ and then iterates $\mathbb{M}_{\mathscr{U}}$ for $\tau$-many times where $\tau=\cf(\tau)<\nu$ then one obtains $\binom{\nu}{\lambda}\rightarrow\binom{\nu}{\lambda}$ in the generic extension (and concomitantly $\binom{\mu}{\lambda}\rightarrow\binom{\mu}{\lambda}$ if $\lambda<\mu=\cf(\mu)<\tau$).
The proof is identical to the proof of the above claim.

We conclude this section with a theorem from \cite{MR2987137} which enables us to have a good control on true cofinalities of certain sequences of regular cardinals.

\begin{theorem}
  \label{thmfrom949} Let $\lambda$ be a Laver-indestructible supercompact cardinal, and assume that $\lambda^{+3}<\nu=\cf(\nu)$ where $\nu^{<\lambda}=\nu$.
  There is a forcing notion $\mathbb{Q}$ which preserves cardinals (but not cofinalities) such that if $G\subseteq\mathbb{Q}$ is $V$-generic then the following statements hold in $V[G]$:
  \begin{enumerate}
    \item [$(a)$] $\lambda>\cf(\lambda)=\omega$.
    \item [$(b)$] $\lambda<\Upsilon_\ell<\nu$ for $\ell\in\{0,1,2\}$.
    \item [$(c)$] $(\lambda_n:n\in\omega)$ is an increasing sequence of regular cardinals and $\lambda=\bigcup_{n\in\omega}\lambda_n$.
    \item [$(d)$] $2^{\lambda_n}=\lambda_n^+$ and $2^{\lambda_n^+}=\lambda_n^{++}$ for every $n\in\omega$.
    \item [$(e)$] $tcf(\prod_{n\in\omega}\lambda_n^{+\ell},J)=\Upsilon_\ell$ for $\ell\in\{0,1,2\}$ for some $J\supseteq J^{bd}_\omega$.
  \end{enumerate}
\end{theorem}
We indicate that the $\Upsilon_\ell$s are regular cardinals, but they need not be distinct, nor satisfy some order relation between them.
We also indicate that $\mathbb{Q}$ is $\lambda^+$-cc.

\newpage 

\section{Terraced cubes and generalized terraced cubes}

For an infinite cardinal $\lambda$ let us call the cube relation $\left( \begin{smallmatrix} \lambda^{++} \\ \lambda^+ \\ \lambda \end{smallmatrix} \right) \rightarrow \left( \begin{smallmatrix} \lambda^{++} \\ \lambda^+ \\ \lambda \end{smallmatrix} \right)$ the \emph{terraced cube relation} at $\lambda$.
Under \textsf{AD} the terraced relation holds at many infinite cardinals as proved in \cite{MR4101445}.
It is unknown whether $\binom{\lambda^{++}}{\lambda^+}\rightarrow\binom{\lambda^{++}}{\lambda^+}$ is consistent with \textsf{ZFC} for some $\lambda$, so we do not have a model in which we can even start checking the possibility of the positive terraced cube relation at $\lambda$.

It has been proved in \cite{MR0202613} that $2^\lambda=\lambda^+$ implies $\binom{\lambda^+}{\lambda}\nrightarrow\binom{\lambda^+}{\lambda}$, so a necessary assumption for the terraced relation is $2^\lambda>\lambda^+$ and $2^{\lambda^+}>\lambda^{++}$.
In this section we show that in some sense these assumptions are far from being sufficient.
We shall prove that if $2^\lambda=\lambda^{++}$ (so the necessary $2^\lambda>\lambda^+$ is satisfied) then $\left( \begin{smallmatrix} \lambda^{++} \\ \lambda^+ \\ \lambda \end{smallmatrix} \right) \nrightarrow \left( \begin{smallmatrix} \lambda^{++} \\ \lambda^+ \\ \lambda \end{smallmatrix} \right)$ no matter how large is $2^{\lambda^+}$ or any other property of $\lambda^+$.

\begin{theorem}
\label{thmterraced} Let $\lambda$ be an infinite cardinal. \newline
If $2^\lambda\leq\lambda^{++}$ then $\left( \begin{smallmatrix} \lambda^{++} \\ \lambda^+ \\ \lambda \end{smallmatrix} \right) \nrightarrow \left( \begin{smallmatrix} \lambda^{++} \\ \lambda^+ \\ \lambda \end{smallmatrix} \right)$.
\end{theorem}

\par\noindent\emph{Proof}. \newline
If $2^\lambda=\lambda^+$ then $\binom{\lambda^+}{\lambda}\nrightarrow\binom{\lambda^+}{\lambda}$ and the failure of the terraced relation follows.
Assume, therefore, that $2^\lambda=\lambda^{++}$ and let $\{C_\zeta:\zeta\in\lambda^{++}\}$ be an enumeration of $[\lambda]^\lambda$.
For every $\alpha\in\lambda^{++}$ let $\mathcal{F}_\alpha=\{C_\zeta:\zeta\in\alpha\}$, so $|\mathcal{F}_\alpha|\leq\lambda^+$.
Let $\{C^\alpha_\zeta:\zeta\in\lambda^+\}$ be an enumeration of $\mathcal{F}_\alpha$ of order type $\lambda^+$, using repetitions if needed.
We may also assume, without loss of generality, that $\bigcup\mathcal{F}_\alpha=\lambda$.

By induction on $\alpha\in\lambda^{++}$ we define $d_\alpha:\{\alpha\}\times\lambda^+\times\lambda\rightarrow 2$.
Arriving at the ordinal $\alpha$ we define for every $\beta\in\lambda^+$ the family $\mathcal{G}^\alpha_\beta=\{C^\alpha_\zeta:\zeta\in\beta\}$.
Since $|\beta|\leq\lambda$ we can reenumerate the elements of $\mathcal{G}^\alpha_\beta$ by $\{C^\alpha_{\beta\eta}:\eta\in\lambda\}$, possibly with repetitions.

Now for every $\beta\in\lambda^+$ we choose by induction on $\eta\in\lambda$ two distinct elements $i^\alpha_{\beta\eta},j^\alpha_{\beta\eta}\in C^\alpha_{\beta\eta}$ so that $i^\alpha_{\beta\eta},j^\alpha_{\beta\eta}\notin \{i^\alpha_{\beta\sigma},j^\alpha_{\beta\sigma}:\sigma\in\eta\}$.
The choice is possible since $|C^\alpha_{\beta\eta}|=\lambda$ and $\eta\in\lambda$.
We define $d_\alpha(\alpha,\beta,i^\alpha_{\beta\eta})=0$ and $d_\alpha(\alpha,\beta,j^\alpha_{\beta\eta})=1$.
At the end of the process, if $d_\alpha(\alpha,\gamma,\delta)$ is not defined yet for some $(\gamma,\delta)\in\lambda^+\times\lambda$ then we let $d_\alpha(\alpha,\gamma,\delta)=0$.

Define $d:\lambda^{++}\times\lambda^+\times\lambda\rightarrow 2$ by $d=\bigcup\{d_\alpha:\alpha\in\lambda^{++}\}$.
Suppose that $A\in[\lambda^{++}]^{\lambda^{++}},B\in[\lambda^+]^{\lambda^+}$ and $C\in[\lambda]^\lambda$.
Choose $\alpha\in A$ such that $C\in\mathcal{F}_\alpha$, so $C=C^\alpha_\zeta$ for some $\zeta\in\lambda^+$.
Choose $\beta\in B$ such that $C^\alpha_\zeta\in\mathcal{G}^\alpha_\beta$, so $C=C^\alpha_\zeta=C^\alpha_{\beta\eta}$ for some $\eta\in\lambda$.
By the construction of $d$ there are $i^\alpha_{\beta\eta},j^\alpha_{\beta\eta}\in C$ so that $d_\alpha(\alpha,\beta,i^\alpha_{\beta\eta})=0$ and $d_\alpha(\alpha,\beta,j^\alpha_{\beta\eta})=1$, hence $d''(A\times B\times C)=\{0,1\}$ as required.

\hfill \qedref{thmterraced}

Notice that $\left( \begin{smallmatrix} \lambda^{++} \\ \lambda^+ \\ \lambda^{++} \end{smallmatrix} \right) \nrightarrow \left( \begin{smallmatrix} \lambda^{++} \\ \lambda^+ \\ \lambda \end{smallmatrix} \right)$ follows from $2^\lambda\leq\lambda^{++}$ by a similar argument.
One has to replace $[\lambda]^\lambda$ in the above proof by $[\lambda^{++}]^\lambda$ whose size is $\lambda^{++}$ under the assumption $2^\lambda\leq\lambda^{++}$.
Another generalization of the above theorem gives a negative $n$-cube relation whenever $2^\lambda<\lambda^{+n}$.
We phrase the following:

\begin{question}
\label{qterraced} Is the positive strong terraced relation consistent with \textsf{ZFC} for some infinite cardinal $\lambda$?
\end{question}

We make the comment that if the negative terraced relation at $\lambda$ follows from the assumption that $2^\lambda=\lambda^{+3}$ then a negative answer to the above problem will be proved.
Indeed, if $2^\lambda>\lambda^{+3}$ and the positive terraced relation holds at $\lambda$ then it will be preserved upon collapsing $2^\lambda$ to $\lambda^{+3}$.

We mentioned the fact that the negative relation $\binom{\lambda^+}{\lambda}\nrightarrow\binom{\lambda^+}{\lambda}$ follows from the assumption $2^\lambda=\lambda^+$, as proved by Erd\H{o}s, Hajnal and Rado in \cite{MR0202613}.
Naturally, they asked whether one can replace $\lambda^+$ by $2^\lambda$ and prove, in \textsf{ZFC}, that $\binom{2^\lambda}{\lambda}\nrightarrow\binom{2^\lambda}{\lambda}$.
It turned out that the answer is negative, and the positive relation $\binom{2^\lambda}{\lambda}\rightarrow\binom{2^\lambda}{\lambda}$ is consistent in many cases, see \cite{MR2927607}.
An example of this positive relation will be seen later in the paper.

In the light of the negative terraced relation of Theorem \ref{thmterraced} one may wonder about cube relations where the small parameter is $\lambda$ and the large parameter $\nu$ is above $\lambda^{++}$.
In the next section we prove a combinatorial theorem which suggests a way to produce a positive relation.
This will be followed by a forcing construction which gives one type of such a relation, where $2^\lambda$ is relatively small. \\

We move now to a generalized form of terraced relations.
Suppose that $\lambda<\mu<\nu\leq 2^\mu$.
A cube relation of the form $\left( \begin{smallmatrix} \nu \\ \mu \\ \lambda \end{smallmatrix} \right) \rightarrow \left( \begin{smallmatrix} \alpha \\ \beta \\ \gamma \end{smallmatrix} \right)$ will be called \emph{strong} if $\alpha=\nu,\beta=\mu$ and $\gamma=\lambda$.
Strong terraced cube relations are the special case in which $\mu=\lambda^+$ and $\nu=\lambda^{++}$.
As in the previous section we focus on strong relations of the form $\left( \begin{smallmatrix} \nu \\ \mu \\ \lambda \end{smallmatrix} \right) \rightarrow \left( \begin{smallmatrix} \nu \\ \mu \\ \lambda \end{smallmatrix} \right)$, namely the size of the monochromatic cube is the size of the domain of the coloring.

We shall prove that under large cardinal assumptions one obtains the consistency of a positive strong cube relation over $\lambda$.
In the light of the result from the previous section, a reasonable assumption to begin with is $\nu>\lambda^{++}$.
This will be discussed at the end of the proof.
The small component $\lambda$ in our theorem will be a strong limit singular cardinal.
In this section we shall prove that under some combinatorial assumptions one has a positive relation.
In the next section we shall see that these assumptions are forceable.

In the theorem below, the cofinality of $\lambda$ is $\omega$.
Countable cofinality is not essential for the main result, and this issue will be discussed briefly at the end of this section.
The proof of the main theorem of this section is a modification of \cite[Theorem 5.1]{MR3509813}.
However, we need a bit more in order to get the cube relation.

\begin{theorem}
\label{thmechelon} Assume that:
\begin{enumerate}
\item [$(a)$] $\lambda>\cf(\lambda)=\omega$ and $\mathsf{CH}$ holds.
\item [$(b)$] $\lambda<\mu=\cf(\mu)<\nu=\cf(\nu)\leq 2^\mu$.
\item [$(c)$] $(\lambda_n:n\in\omega)$ is an increasing sequence of regular cardinals so that $\lambda=\bigcup_{n\in\omega}\lambda_n$, and $\mathsf{GCH}$ holds at $\lambda_n$ and $\lambda_n^+$ for every $n\in\omega$.
\item [$(d)$] $\Upsilon_\ell={\rm tcf}(\prod_{n\in\omega}\lambda_n^{+\ell},J^{\rm bd}_\omega)$ for $\ell\in\{0,1,2\}$.
\item [$(e)$] $\mu<\Upsilon_\ell$ and $\nu>\Upsilon_\ell$ for $\ell\in\{0,1,2\}$.
\item [$(f)$] $\binom{\nu}{\mu}\rightarrow\binom{\nu}{\mu}_{\aleph_1}$.
\end{enumerate}
Then $\left( \begin{smallmatrix} \nu \\ \mu \\ \lambda \end{smallmatrix} \right) \rightarrow \left( \begin{smallmatrix} \nu \\ \mu \\ \lambda \end{smallmatrix} \right)$.
\end{theorem}

\par\noindent\emph{Proof}. \newline
Suppose that $d:\nu\times\mu\times\lambda\rightarrow 2$ is given.
We fix three scales as follows.
Let $\bar{f}=(f_\alpha:\alpha\in\Upsilon_2)$ be a scale in $\prod_{n\in\omega}\lambda_n^{++}$, let $\bar{g}=(g_\varepsilon:\varepsilon\in\Upsilon_1)$ be a scale in $\prod_{n\in\omega}\lambda_n^+$ and let $\bar{h}=(h_\delta:\delta\in\Upsilon_0)$ be a scale in $\prod_{n\in\omega}\lambda_n$.
For every $\alpha\in\nu,\beta\in\mu,n\in\omega$ and $i\in\{0,1\}$, define:
$$
A^i_{\alpha\beta n}=\{\gamma\in\lambda_n^+:d(\alpha,\beta,\gamma)=i\}.
$$
For every $n\in\omega$ let $\{S^n_j:j\in\lambda_n^{++}\}$ be an enumeration of $\mathcal{P}(\lambda_n^+)$.
For every $\alpha\in\nu,\beta\in\mu,i\in\{0,1\}$ we define a function $f^i_{\alpha\beta}\in\prod_{n\in\omega}\lambda_n^{++}$ by $f^i_{\alpha\beta}(n)=\min\{j\in\lambda_n^{++}:A^i_{\alpha\beta n}=S^n_j\}$.
The color $i$ can be removed by defining $f_{\alpha\beta}(n)=\max\{f^i_{\alpha\beta}(n):i<2\}$.

Thus we have defined $\mu$ functions for each $\alpha\in\nu$ of the form $f_{\alpha\beta}$.
Applying assumption $(e)$ we can choose $f_\alpha\in\bar{f}$ so that $\beta\in\mu\Rightarrow f_{\alpha\beta}<_{J^{\rm bd}_\omega}f_\alpha$, and we do this for every $\alpha\in\nu$.
Since $\Upsilon_2<\nu=\cf(\nu)$, there exists $A_0\in[\nu]^\nu$ and a fixed $f\in\bar{f}$ so that $\alpha\in A_0\Rightarrow f_\alpha=f$.
By increasing each $f(n)$ if needed we may assume that $f(n)\geq\lambda_n^+$ for every $n\in\omega$.

For every $\alpha\in{A_0},\beta\in\mu$ we choose $n_{\alpha\beta}\in\omega$ so that $n\geq n_{\alpha\beta}\Rightarrow f_{\alpha\beta}(n)<f(n)$.
The mapping $(\alpha,\beta)\mapsto n_{\alpha\beta}$ admits a monochromatic product of size $\nu\times\mu$ by virtue of $(f)$, so choose $A_1\in[A_0]^\nu,B_1\in[\mu]^\mu$ and $n_1\in\omega$ such that $(\alpha,\beta)\in A_1\times B_1\Rightarrow n_{\alpha\beta}=n_1$.

Observe that $|\{S^n_j:j\in f(n)|=\lambda_n^+$ for every $n\in\omega$, so we can reenumerate these sets by $\{T^n_j:j\in\lambda_n^+\}$.
For every $\alpha\in A_1,\beta\in B_1$ and $i\in\{0,1\}$ define $g^i_{\alpha\beta}\in\prod_{n\in\omega}\lambda_n^+$ by $g^i_{\alpha\beta}(n)=\min\{j\in\lambda_n^+:A^i_{\alpha\beta n}=T^n_j\}$ whenever $n\geq{n_1}$ and $g^i_{\alpha\beta}(n)=0$ when $n<n_1$.
Let $g_{\alpha\beta}(n)=\max\{g^i_{\alpha\beta}(n):i<2\}$ for every $n\in\omega$.
By the same reasoning as before we choose $A_2\in[A_1]^\nu,B_2\in[B_1]^\mu,n_2\in\omega$ and $g\in\bar{g}$ such that $\alpha\in A_2,\beta\in B_2\Rightarrow g_{\alpha\beta}<_{J^{\rm bd}_\omega}g$. Moreover, $g(n)\geq\lambda_n$ for every $n\in\omega, n_1\leq n_2$ and $g_{\alpha\beta}(n)<g(n)$ whenever $\alpha\in A_2,\beta\in B_2$ and $n\geq n_2$.
These statements follow from the properties of $\bar{g}$ and assumption $(f)$.

We need another round of the same process, so we focus now on the set $\{T^n_j:j\in g(n)\}$ for every $n\in\omega$ whose size is $\lambda_n$.
By reenumerating these sets as $\{W^n_j:j\in\lambda_n\}$ we define for each $\alpha\in A_2,\beta\in B_2$ and $i\in\{0,1\}$ the function $h^i_{\alpha\beta}\in\prod_{n\in\omega}\lambda_n$ in a similar fashion.
Namely, if $n\geq{n_2}$ then $h^i_{\alpha\beta}(n)=\min\{j\in\lambda_n:A^i_{\alpha\beta n}=W^n_j\}$, if $n<n_2$ then $h^i_{\alpha\beta}(n)=0$ and then $h_{\alpha\beta}(n)=\max\{h^i_{\alpha\beta}(n):i<2\}$ for every $n\in\omega$.
At last, we shrink ourselves to $A_3\in[A_2]^\nu,B_3\in[B_2]^\mu$ and we choose $n_2\leq n_3\in\omega$ and $h\in\prod_{n\in\omega}\lambda_n,h\in\bar{h}$ so that $\alpha\in A_3,\beta\in B_3\Rightarrow h_{\alpha\beta}<_{J^{\rm bd}_\omega}h$ and $h_{\alpha\beta}(n)<h(n)$ whenever $\alpha\in A_3,\beta\in B_3$ and $n\geq n_3$.

For every $n\in\omega$ we define an equivalence relation $e_n$ on the ordinals of $\lambda_n^+$ as follows:
$$
\gamma_0 e_n \gamma_1\quad \text{iff}\quad \forall j\in h(n), \gamma_0\in W^n_j\Leftrightarrow\gamma_1\in W^n_j.
$$
Since we are assuming \textsf{GCH} at the $\lambda_n$s and their successors, we see that the number of equivalence classes of $e_n$ is less than $\lambda_n^+$.
Hence one can choose for every $n\in\omega$ an equivalence class $E_n$ of $e_n$ and a color $i^n_{\alpha\beta}\in\{0,1\}$ for every $\alpha\in A_3,\beta\in B_3$ such that $|E_n|=\lambda_n^+$ and $\gamma\in E_n\Rightarrow d(\alpha,\beta,\gamma)=i^n_{\alpha\beta}$.

For every $\alpha\in A_3,\beta\in B_3$ we choose an infinite set $u_{\alpha\beta}\subseteq\omega$ and a color $i_{\alpha\beta}\in\{0,1\}$ such that $n\in u_{\alpha\beta}\Rightarrow i^n_{\alpha\beta}=i_{\alpha\beta}$.
Applying assumption $(f)$ we can find $A\in[A_3]^\nu,B\in[B_3]^\mu,u\in[\omega]^\omega$ and a fixed color $i\in\{0,1\}$ such that $\alpha\in A,\beta\in B,n\in u\Rightarrow i_{\alpha\beta}=i$.
We indicate that here we use assumption $(a)$ and the fact that we have $\aleph_1$-many colors in assumption $(f)$.
Let $C=\bigcup\{E_n:n_3\leq n\in u\}$, so $|C|=\mu$ as $u$ is unbounded in $\omega$.
It follows that $d''(A\times B\times C)=\{i\}$, so the proof is accomplished.

\hfill \qedref{thmechelon}

\begin{remark}
  \label{rafter23}
  The assumption $\mu<\Upsilon_\ell$ for $\ell\in\{0,1,2\}$ in part $(e)$ of the above theorem can be replaced by the weaker assumption $\mu\neq\Upsilon_\ell$.
  To see this, consider the choice of $f_\alpha$ for every $\alpha\in\nu$, based on $(e)$.
  In the above proof we choose $f_\alpha\in\bar{f}$ so that $\beta\in\mu\Rightarrow f_{\alpha\beta}<_{J^{bd}_\omega}f_\alpha$.
  However, it is sufficient to find $C\in[\mu]^\mu$ so that $\beta\in{C}\Rightarrow f_{\alpha\beta}<_{J^{bd}_\omega}f_\alpha$ and then to focus on the elements of $C$ from this point onwards.
  The same holds at the choice of the functions $g$ and $h$.
\end{remark}

\hfill \qedref{rafter23}

The above theorem is based on several assumptions, and one may wonder whether these assumptions are forceable.
In the next section we shall see that if there are two supercompact cardinals in the ground model then the answer is positive.

We make the comment that one can modify the proof and incorporate singular cardinals with uncountable cofinality.
The required changes are choosing an appropriate increasing sequence $(\lambda_\varepsilon:\varepsilon\in\cf(\lambda))$, requiring \textsf{GCH} at $\lambda_\varepsilon$ and $\lambda_\varepsilon^+$ for every $\varepsilon\in\cf(\lambda)$, and increasing the number of colors in assumption $(f)$ from $\aleph_1$ to $2^{\cf(\lambda)}$.

\newpage 

\section{Forcing our assumptions}

As indicated in the introduction, in order to prove the consistency of $\left( \begin{smallmatrix} \nu \\ \mu \\ \lambda \end{smallmatrix} \right) \rightarrow \left( \begin{smallmatrix} \nu \\ \mu \\ \lambda \end{smallmatrix} \right)$ one has to make sure that all possible pairs satisfy the pertinent positive relation.
The relations in which $\lambda$ is involved are relatively easy to force when $\lambda$ is a strong limit singular cardinal.
This is done in \cite[Theorem 5.1]{MR3509813}.
Actually, the proof of the main theorem in the previous section is based on the proof of that theorem, upon adding the assumption $\binom{\nu}{\mu}\rightarrow\binom{\nu}{\mu}_{\aleph_1}$ and handling cubes.
So basically one has to force the assumptions of \cite[Theorem 5.1]{MR3509813} and concomitantly the additional assumption $(f)$.

The challenging relation which we need is the relation $\binom{\nu}{\mu}\rightarrow\binom{\nu}{\mu}$, and what is more $\binom{\nu}{\mu}\rightarrow\binom{\nu}{\mu}_{\aleph_1}$.
If $\mu$ is a singular cardinal with small cofinality then $\binom{\nu}{\mu}\rightarrow\binom{\nu}{\mu}$ will follow from $\binom{\nu}{\cf(\mu)}\rightarrow\binom{\nu}{\cf(\mu)}$ which can be easily arranged.
Indeed, one can choose $\mu$ so that $\cf(\mu)=\omega$ and force $2^{\aleph_0}<\lambda$.
In such a case, the relation $\binom{\nu}{\cf(\mu)}\rightarrow\binom{\nu}{\cf(\mu)}$ follows trivially since $2^{\aleph_0}<\nu$.
However, this is not what the poet meant.
Hence we are asking for a regular cardinal $\mu$ (and similarly, a regular cardinal $\nu$) with respect to our cube relation.
In such cases, it is harder to force $\binom{\nu}{\mu}\rightarrow\binom{\nu}{\mu}$.

It has been proved in \cite{MR3201820} that if $\mathfrak{r}_\mu<\nu$ and $\cf(\nu)>\mathfrak{r}_\mu$ then $\binom{\nu}{\mu}\rightarrow\binom{\nu}{\mu}$.
However, we need more than two colors which makes life a bit more complicated.
Rather than $\mathfrak{r}_\mu$ we shall work with $\mathfrak{u}_\mu$ (these characteristics are similar and actually it is unknown whether they can be separated where $\mu>\aleph_0$).
Recall that a base of a uniform ultrafilter $\mathscr{U}$ over $\mu$ is a subset $\mathcal{B}$ of $\mathscr{U}$ such that for every $A\in\mathscr{U}$ one can find $B\in\mathcal{B}$ for which $B\subseteq A$.
The character ${\rm Ch}(\mathscr{U})$ is the minimal cardinality of a base of $\mathscr{U}$.
The ultrafilter number $\mathfrak{u}_\mu$ is the minimal value of ${\rm Ch}(\mathscr{U})$ for some uniform ultrafilter over $\mu$.
If we restrict our attention to complete ultrafilters then we can get more colors from the assumption $\mathfrak{u}_\mu<2^\mu$ as mirrored by the following:

\begin{claim}
\label{clmufnumber} Suppose that $\mathscr{U}$ is a $\theta$-complete uniform ultrafilter over $\mu, \partial={\rm Ch}(\mathscr{U})<\nu=\cf(\nu)$ and $\chi<\theta$. Then $\binom{\nu}{\mu}\rightarrow\binom{\nu}{\mu}_\chi$.
\end{claim}

\par\noindent\emph{Proof}. \newline
Let $\mathcal{B}=\{B_\delta:\delta\in\partial\}\subseteq\mathscr{U}$ be a base for $\mathscr{U}$.
Suppose that $c:\nu\times\mu\rightarrow\chi$ is a coloring.
For every $\alpha\in\nu$ and every $i\in\chi$ let $A^i_\alpha=\{\beta\in\mu:c(\alpha,\beta)=i\}$.
Since $\mathscr{U}$ is $\theta$-complete and $\chi<\theta$ one can find for every $\alpha\in\nu$ a color $i(\alpha)$ such that $A^{i(\alpha)}_\alpha\in\mathscr{U}$.
Since $\mathcal{B}$ is a base of $\mathscr{U}$ one can find $\delta(\alpha)\in\partial$ such that $B_{\delta(\alpha)}\subseteq A^{i(\alpha)}_\alpha$.

Choose $A'\in[\nu]^\nu$ and a fixed color $i\in\chi$ such that $\alpha\in A'\Rightarrow i(\alpha)=i$.
Since $\partial<\nu=\cf(\nu)$ one can find a fixed $\delta\in\partial$ and $A\in[A']^\nu$ such that $\alpha\in A\Rightarrow\delta(\alpha)=\delta$.
Let $B=B_\delta$, and recall that $|B|=\mu$ since $B\in\mathscr{U}$.
But now we are done since $c''(A\times B)=\{i\}$.

\hfill \qedref{clmufnumber}

The ability to force $\mathfrak{u}_\mu<2^\mu$ and moreover with ultrafilters which possess some degree of completeness is supplied by \cite{MR4081063}.
Suppose that $\kappa$ is a Laver-indestructible supercompact cardinal, $\lambda'\leq\kappa$ and $\mu>\kappa$.
Assume further that $\mathscr{D}$ is a $\lambda'$-complete $\theta$-indecomposable uniform filter over $\mu$, where $\theta=\cf(\mu')$ and $\mu'$ is a strong limit singular cardinal above $\mu$.
It is shown in \cite{MR4081063} that there are forcing notions which make $2^\mu>\mu'$ and every uniform ultrafilter $\mathscr{U}$ which extends $\mathscr{D}$ in the generic extension satisfies ${\rm Ch}(\mathscr{U})\leq\mu'<2^\mu$.

Now if $\kappa$ is Laver indestructible and the forcing $\mathbb{P}$ is $\kappa$-directed-closed then $\kappa$ remains supercompact in $V[G]$ where $G\subseteq\mathbb{P}$ is $V$-generic.
Therefore, one can choose a uniform $\lambda'$-complete ultrafilter which extends $\mathscr{D}$ and apply Claim \ref{clmufnumber}.
Our strategy will be to force $\binom{\nu}{\mu}\rightarrow\binom{\nu}{\mu}_{\lambda'}$ for some $\lambda'$ and then extend the universe once again in order to obtain the appropriate pcf structure.
It is important at this stage to make sure that the relation $\binom{\nu}{\mu}\rightarrow\binom{\nu}{\mu}_{\lambda'}$ is preserved by the second step of our forcing.
This will be ensured by the following:

\begin{lemma}
\label{lemccc} Assume that $\binom{\nu}{\mu}\rightarrow\binom{\nu}{\mu}_{\lambda'}$ where $\lambda'=\lambda'^{<\lambda'}$.
Let $\mathbb{P}$ be a $\lambda$-cc forcing notion where $\lambda<\lambda'$ and assume further that $|\mathbb{P}|^{<\lambda}\leq\lambda'$.
Let $G\subseteq\mathbb{P}$ be generic over $V$.
Then the positive relation $\binom{\nu}{\mu}\rightarrow\binom{\nu}{\mu}_{\lambda'}$ holds in $V[G]$.
\end{lemma}

\par\noindent\emph{Proof}. \newline
Let $c:\nu\times\mu\rightarrow\lambda'$ be a new coloring and let $\name{c}$ be a $\mathbb{P}$-name for $c$.
For every $\alpha\in\nu,\beta\in\mu$ let $\mathcal{A}_{\alpha\beta}$ be a maximal antichain of conditions which force a value to $\name{c}(\alpha,\beta)$.
Define $f:\nu\times\mu\rightarrow[\mathbb{P}]^{<\lambda}$ by $f(\alpha,\beta)=\mathcal{A}_{\alpha\beta}$, so $f\in V$.

Since $|\mathbb{P}|^{<\lambda}\leq\lambda'$ and $\binom{\nu}{\mu}\rightarrow\binom{\nu}{\mu}_{\lambda'}$ in $V$, there are $A_0\in[\nu]^\nu,B_0\in[\mu]^\mu$ and a maximal antichain $\mathcal{A}\in[\mathbb{P}]^{<\lambda}$ such that $f''(A_0\times B_0)=\mathcal{A}$.
Let $p$ be the unique condition in $G\cap\mathcal{A}$.
Define a coloring $d:A_0\times{B_0}\rightarrow\lambda'$ in $V$, by letting $d(\alpha,\beta)=\xi$ iff $p\Vdash\name{c}(\alpha,\beta)=\xi$.
Apply the relation $\binom{\nu}{\mu}\rightarrow\binom{\nu}{\mu}_{\lambda'}$ in the ground model to obtain $A\in[A_0]^\nu$ and $B\in[B_0]^\mu$ so that $d\upharpoonright(A\times{B})$ is constant.
Since $p\Vdash\name{c}\upharpoonright(A\times{B})=d\upharpoonright(A\times{B})$, we are done.

\hfill \qedref{lemccc}

Equipped with the above lemma, we can phrase and prove the main result of this section:

\begin{theorem}
\label{thmforcemt} Assuming the existence of two supercompact cardinals in the ground model, one can force the cube relation $\left( \begin{smallmatrix} \nu \\ \mu \\ \lambda \end{smallmatrix} \right) \rightarrow \left( \begin{smallmatrix} \nu \\ \mu \\ \lambda \end{smallmatrix} \right)$ when $\lambda<\mu=\cf(\mu)<\nu=\cf(\nu)=2^\mu$.
\end{theorem}

\par\noindent\emph{Proof}. \newline
We commence with a pair of supercompact cardinals $\lambda<\kappa$, where both are Laver-indestructible, according to \cite{MR0472529}.
We fix a regular cardinal $\mu>\kappa^+$, a strong limit singular cardinal $\mu'>\mu$ such that $\kappa<\theta=\cf(\mu')<\mu$ and such that there exists a $\kappa$-complete $\theta$-indecomposable filter $\mathscr{D}$ over $\mu$.
Finally, we fix a regular cardinal $\nu>\mu'$.

Let $\mathbb{P}=\mathcal{C}_{\theta\nu}$ be the usual Cohen forcing that adds $\nu$-many Cohen subsets of $\theta$.
Let $G\subseteq\mathbb{P}$ be $V$-generic.
In $V[G]$, let $\mathscr{U}$ be a $\kappa$-complete ultrafilter over $\mu$ that extends $\mathscr{D}$, we use the fact that $\kappa$ is supercompact in $V[G]$ since $\mathbb{P}$ is $\kappa$-directed-closed.

From \cite{MR4081063} we deduce that $\mathfrak{u}_\mu<2^\mu=\nu$ in $V[G]$.
The positive relation $\binom{\nu}{\mu}\rightarrow\binom{\nu}{\mu}_{<\kappa}$ in $V[G]$ follows from Claim \ref{clmufnumber}.
We define in $V[G]$ a forcing notion $\mathbb{Q}$ of cardinality $\kappa$.
The forcing notion $\mathbb{Q}$ is an iteration in which the first step forces $2^\lambda=\chi$ for some (arbitrarily large) $\chi<\kappa$, and the rest of the iteration is based on \cite{MR2987137}, see Theorem \ref{thmfrom949}.

If $H\subseteq\mathbb{Q}$ is $V[G]$-generic then $2^\lambda=\chi$ and $\lambda$ becomes a strong limit singular cardinal of countable cofinality in $V[G][H]$.
Likewise, $\Upsilon_\ell\neq\mu$ for $\ell\in\{0,1,2\}$, where the $\Upsilon_\ell$s are defined in the assumptions of Theorem \ref{thmechelon}.
Notice that $\mathbb{Q}$ is $\lambda^+$-cc, and $|\mathbb{Q}|^{<\lambda^+}=\lambda'<\kappa$ in $V[G]$.
Thus Lemma \ref{lemccc} applies, where $\mathbb{Q}$ here stands for $\mathbb{P}$ there.
It follows that $\binom{\nu}{\mu}\rightarrow\binom{\nu}{\mu}_{\aleph_1}$ holds in $V[G][H]$, hence the positive cube relation $\left( \begin{smallmatrix} \nu \\ \mu \\ \lambda \end{smallmatrix} \right) \rightarrow \left( \begin{smallmatrix} \nu \\ \mu \\ \lambda \end{smallmatrix} \right)$ holds there as well by the conclusion of Theorem \ref{thmechelon}.

\hfill \qedref{thmforcemt}

Let us indicate that the most interesting case in this setting is $\mu=\kappa^{+3}$ and $\theta=\kappa^{++}$.
It is still tempting to adjust the above results to the case in which $2^\lambda=\nu$.
We conclude, therefore, with a couple of open problems.

\begin{question}
\label{qbasa} Is it consistent that $\lambda<\mu=\cf(\mu)<\nu=\cf(\nu)\leq 2^\lambda$ and $\left( \begin{smallmatrix} \nu \\ \mu \\ \lambda \end{smallmatrix} \right) \rightarrow \left( \begin{smallmatrix} \nu \\ \mu \\ \lambda \end{smallmatrix} \right)$ holds?
\end{question}

It seems that if $\lambda=\aleph_0$ then it should be easier to force a positive relation, hence we raise the following:

\begin{question}
\label{qctble} Is it consistent that $\left( \begin{smallmatrix} \nu \\ \mu \\ \omega \end{smallmatrix} \right) \rightarrow \left( \begin{smallmatrix} \nu \\ \mu \\ \omega \end{smallmatrix} \right)$ where $\mu,\nu$ are regular and $\nu\leq 2^\omega$?
\end{question}

\newpage

\bibliographystyle{alpha}
\bibliography{arlist}

\end{document}